\documentclass[oneside]{amsart}
\usepackage[mathscr]{eucal}
\usepackage{amssymb, amsmath, array, amscd}
\usepackage[dvips]{graphicx}
\newtheorem{theorem}{Theorem}[section]

\newtheorem{proposition}{Proposition}[section]
\newtheorem{lemma}{Lemma}[section]

\newtheorem{remark}{Remark}[section]
\newtheorem{example}{Example}[section]
\newtheorem{conjecture}{Conjecture}[section]

\begin{document}

\title[DYNAMICS OF CODIMENSION ONE FOLIATIONS]
{ON THE DYNAMICS OF CODIMENSION ONE HOLOMORPHIC FOLIATIONS WITH
AMPLE NORMAL BUNDLE}
\author{Marco Brunella}
\address{Marco Brunella, IMB - CNRS UMR 5584, 9 Avenue Savary,
21078 Dijon, France}
\begin{abstract}
We investigate the accumulation to singular points of leaves of
codimension one foliations whose normal bundle is ample, with
emphasis on the nonexistence of Levi-flat hypersurfaces.
\end{abstract}

\maketitle

\section{Introduction}

Let $X$ be a compact complex manifold and let ${\mathcal F}$ be a
(singular) codimension one holomorphic foliation on $X$. From the
dynamical point of view, one of the most basic problems we may
address is the following \cite{CLS}: does it hold that every leaf of
${\mathcal F}$ accumulates to the singular set $Sing({\mathcal F})$?
Of course, for a positive answer we need some hypothesis on $X$
and/or ${\mathcal F}$. In \cite{Lin}, Lins Neto showed that the
answer to the above question is positive when $X$ is the complex
projective space ${\mathbb C}P^n$ and $n\ge 3$ (the case $n=2$ is
still open, and much more difficult). Our aim is to prove a similar
result on any manifold $X$, under some assumption on ${\mathcal F}$.

Recall that we may associate to ${\mathcal F}$ its {\it normal
bundle} $N_{\mathcal F}$: the foliation is locally defined by
integrable holomorphic 1-forms $\omega_j\in\Omega^1(U_j)$, with zero
set of codimension at least two, and so $N_{\mathcal F}$ is the
holomorphic line bundle defined by the cocycle $\{ g_{jk}\in
{\mathcal O}^*(U_j\cap U_k)\}$ given by $\omega_j = g_{jk}\omega_k$.
Our feeling is that $N_{\mathcal F}$ should reflect some dynamical
properties of ${\mathcal F}$. A first instance of this philosophy is
the result of \cite{BLM} concerning the existence of loops with
hyperbolic holonomy: by an easy ``regularisation'' of their proof,
one can see that such a result holds under the sole assumption that
$N_{\mathcal F}$ is {\it ample}. Other results involving some form
of positivity of $N_{\mathcal F}$ can be found in \cite{Der} and
\cite{Bru}.

Let us put forward a conjecture.

\begin{conjecture}\label{conj1}
Let $X$ be a compact connected complex manifold of dimension $n\ge
3$, and let ${\mathcal F}$ be a codimension one holomorphic
foliation on $X$ whose normal bundle $N_{\mathcal F}$ is ample. Then
every leaf of ${\mathcal F}$ accumulates to $Sing({\mathcal F})$.
\end{conjecture}

Note that by Baum-Bott formula \cite{Suw} the ampleness of
$N_{\mathcal F}$ implies, at least, that $Sing({\mathcal F})$ is not
empty; we shall return later on this point.

If $X={\mathbb C}P^n$, or more generally if $X$ admits an hermitian
metric of positive curvature, then the ampleness hypothesis is
automatically satisfied, because $N_{\mathcal F}$ is, outside the
singular set, a quotient of $TX$ and therefore it is more positive
than $TX$. Hence the above Conjecture would extend the result of
\cite{Lin}. However, there are many cases in which we can guarantee
that $N_{\mathcal F}$ is ample even if $TX$ is far from positive,
for exemple when $X$ is a complete intersection (of dimension at
least 3) in ${\mathbb C}P^N$ (indeed, $Pic(X)={\mathbb Z}$ by
Lefschetz and thus either $N_{\mathcal F}$ is ample or $N_{\mathcal
F}^*$ is effective, but the latter is excluded by $\Omega^1(X)=\{
0\}$).

Unfortunately, we are unable to prove the above Conjecture in full
generality. Suppose that Conjecture \ref{conj1} does not hold: then
there exists a nonempty compact subset ${\mathcal M}\subset X$ which
is invariant by ${\mathcal F}$ and disjoint from $Sing({\mathcal
F})$. We shall prove that such a ${\mathcal M}$ cannot be a
sufficiently smooth real hypersurface.

More precisely, and more generally:

\begin{theorem}\label{thm1}
Let $X$ be a compact connected K\"ahler manifold of dimension $n\ge
3$. Let $M\subset X$ be a closed real hypersurface of class
$C^{2,\alpha}$, $\alpha > 0$. Suppose that on some neighbourhood $U$
of $M$ there exists a codimension one holomorphic foliation
${\mathcal F}$ which leaves invariant $M$. Then the normal bundle
$N_{\mathcal F}$ does not admit, on $U$, any hermitian metric with
positive curvature.
\end{theorem}

It is perhaps worth noting that, in the situation of this Theorem,
the foliation ${\mathcal F}$ is necessarily nonsingular around $M$.

We stated the result in a slightly cumbersome way because we hope
that it could be generalised to a larger context, by ``shrinking $U$
to $M$'', in the same way as \cite{Siu} is a generalisation of
\cite{Lin}. That is, suppose that $M\subset X$ is a sufficiently
smooth Levi-flat hypersurface (i.e., a real hypersurface smoothly
foliated by complex hypersurfaces), but not necessarily invariant by
a holomorphic foliation on a neighbourhood of it. It is still
possible to define a normal bundle $N_{{\mathcal F}_M}$ of the Levi
foliation ${\mathcal F}_M$ on $M$. It is a smooth ${\mathbb
C}$-bundle on $M$, holomorphic along the leaves; generally speaking,
$N_{{\mathcal F}_M}$ cannot be extended to a holomorphic line bundle
on a neighbourhood of $M$, in the same way as ${\mathcal F}_M$
cannot be extended to a holomorphic foliation \cite{BdB}. However,
it makes sense to say that $N_{{\mathcal F}_M}$ admits a hermitian
metric of positive curvature along the leaves of ${\mathcal F}_M$.
We think that it is {\it never} the case.

When $X={\mathbb C}P^n$ such a generalisation has been done by Siu
in \cite{Siu}; in that case $N_{{\mathcal F}_M}$ {\it does} admit a
metric with leafwise positive curvature, by quotient of the
Fubini-Study metric, and so the conclusion is that $M$ does not
exist.

Remark that if $M$ is real analytic then ${\mathcal F}_M$ and
$N_{{\mathcal F}_M}$ can always be holomorphically extended to some
neighbourhood $U$ of $M$, and a metric on $N_{{\mathcal F}_M}$ with
leafwise positive curvature can always be extended to a metric on
the extended bundle with positive curvature: the positivity in the
direction transverse to the leaves can be gained by multiplying any
extended metric by the factor $\exp [-C dist(\cdot ,M)^2]$, $C\gg
0$. Thus we return to the setting of Theorem \ref{thm1}.

If $M$ is of class $C^\infty$, ${\mathcal F}_M$ can be extended at
least as a ``formal'' object, i.e. as a $C^\infty$-section of the
bundle of hyperplanes of $X$ whose $\bar\partial$ vanishes along $M$
at infinite order. Optimistically, one could try to construct a
holomorphic extension using a sort of formal-convergent principle
{\it \`a la} Hironaka-Matsumura, and the leafwise positivity of
$N_{{\mathcal F}_M}$. This is related to the vanishing theorem
proved in \cite{Siu} and \cite{Bri}, and it remains to see if that
vanishing theorem can be applied in our context.

Similar problems can be posed also in the context of Conjecture
\ref{conj1}. That is, one may ask about the existence of codimension
one laminations (not necessarily invariant by a holomorphic
foliation) whose normal bundle is leafwise ample. Even in the case
of $X={\mathbb C}P^n$, this seems a still open problem.

In another direction, we hope that Theorem \ref{thm1} could be
useful to classify Levi-flat hypersurfaces in Fano manifolds, at
least under the assumption of invariance by a global holomorphic
foliation. Indeed, if $X$ is a Fano manifold (i.e., its
anticanonical bundle $K_X^*$ is ample), then $N_{\mathcal F} =
K_X^*\otimes K_{\mathcal F}$ is more positive than the canonical
bundle $K_{\mathcal F}$ of the foliation. Thus, if $M\subset X$ is a
real hypersurface of class $C^{2,\alpha}$ invariant by ${\mathcal
F}$, then $K_{\mathcal F}$ cannot be nef, by Theorem \ref{thm1}
(recall that nef + ample is ample). According to Miyaoka and
Shepherd-Barron \cite{ShB}, this gives remarkable informations on
${\mathcal F}$, in particular concerning the existence of rational
curves inside the leaves. Eventually, all of this could prove that
$M$ is smoothly fibered by compact rationally connected submanifolds
of dimension $n-2$, contained in the leaves. In a different context,
the one of complex tori, a similar structure has been found by
Ohsawa \cite{Ohs}.

Our proof of Theorem \ref{thm1} follows the same path as \cite{Lin}.
We proceed by contradiction, by assuming that $N_{\mathcal F}$ has a
metric of positive curvature. The main step consists in proving that
the complement $X\setminus M$ is then {\it strongly pseudoconvex},
that is a point modification of a Stein space. This step is by free
in \cite{Lin}, in an even stronger form, thanks to the classical
solution of the Levi problem in projective spaces \cite{Fuj}
\cite{Tak}. Of course, Levi problem has a negative answer on most
compact K\"ahler manifolds, and in our case we need a delicate
glueing procedure of plurisubharmonic functions, which involves the
$C^{2,\alpha}$-regularity of $M$. Then we conclude the proof in two
slightly different ways, one close to \cite{Siu} and the other close
to \cite{Lin}. It is only in this last part that the dimensional
hypothesis $n\ge 3$ is used.

In the last Section we discuss some possible ways to fill the gap
between Theorem \ref{thm1} and Conjecture \ref{conj1}.

\section{Convexity of the complement}

Recall that a complex manifold $V$ is {\it strongly pseudoconvex}
(or {\it 1-convex}) if there exists a $C^2$ function $\psi :V\to
{\mathbb R}$ which is:
\begin{enumerate}
\item[(i)] exhaustive, i.e. $\{ \psi\le c\}$ is compact for every
$c\in{\mathbb R}$;
\item[(ii)] strictly plurisubharmonic, i.e.
${\rm i}\partial\bar\partial\psi > 0$, outside a compact subset.
\end{enumerate}
Classical results by Grauert and Remmert \cite[\S 2]{Pet} say that a
strongly pseudoconvex manifold $V$ is a point modification of a
Stein space $V_0$: there exists a proper holomorphic map $\pi :V\to
V_0$ and a finite subset $P\subset V_0$ such that $\pi$ is an
isomorphism between $V\setminus \pi^{-1}(P)$ and $V_0\setminus P$.
The {\it exceptional subset} $\pi^{-1}(P)$ is the maximal compact
analytic subset of $V$ of positive dimension.

Before starting the proof of Theorem \ref{thm1}, let us recall the
well-known and easy proof of the following model case: if $X$ is a
compact connected complex manifold and $Y\subset X$ is a (smooth)
complex hypersurface whose normal bundle $N_Y$ is ample (on $Y$),
then $X\setminus Y$ is strongly pseudoconvex. Indeed, by the
adjunction formula $N_Y$ may be identified with ${\mathcal
O}_X(Y)\vert_Y$, and so we may construct on ${\mathcal O}_X(Y)$ a
hermitian metric whose curvature is positive on some neighbourhood
of $Y$. The line bundle ${\mathcal O}_X(Y)$ has a global holomorphic
section $s$ on $X$, which vanishes exactly on $Y$. Then the function
$\psi = -\log\| s\|$, where $\|\cdot \|$ is the above metric, is an
exhaustion of $X\setminus Y$ and is strictly plurisubharmonic close
to $Y$. Here $\| s\|$ plays the role of distance from $Y$, and below
we shall need to find a good substitute for it in our context. See
also \cite{Tak} and \cite{Ohs} (without mentioning Oka...).

Consider $X$, $M$, ${\mathcal F}$ as in Theorem \ref{thm1}, and
suppose by contradiction that $N_{\mathcal F}$ has, on some
neighbourhood $U$ of $M$, a hermitian metric with positive curvature
$\omega$. In this Section we shall prove:

\begin{proposition}\label{prop1}
$X\setminus M$ is strongly pseudoconvex.
\end{proposition}

We may assume that $M$ is connected and, up to taking a double
covering, orientable; hence, up to reducing $U$, $U\setminus M$ has
two connected components $U^+$ and $U^-$. We may choose $U$ so that
it is covered by a finite number of charts $U_j\simeq {\mathbb
D}\times {\mathbb D}^{n-1}$ adapted to the foliation, $j=1,\ldots
,\ell$, and $M$ cuts each $U_j$ along $M_j = \gamma_j\times{\mathbb
D}^{n-1}$, where $\gamma_j\subset{\mathbb D}$ is a proper arc of
class $C^{2,\alpha}$. Thus $U_j\setminus M_j$ has two connected
components $U_j^+$ and $U_j^-$, contained respectively in $U^+$ and
$U^-$. We have $U_j^+=V_j^+\times{\mathbb D}^{n-1}$ and
$U_j^-=V_j^-\times{\mathbb D}^{n-1}$, where $V_j^+$ and $V_j^-$ are
the two connected components of ${\mathbb D}\setminus\gamma_j$.

For each $j=1,\ldots ,\ell$, choose a biholomorphism
$$\varphi_j : V_j^+ \longrightarrow {\mathbb D}^+ = \{ z\in{\mathbb
D}\ \vert\ \Im m z > 0\}$$ sending the arc $\gamma_j\subset\partial
V_j^+$ to the arc $(-1,1)\subset\partial{\mathbb D}^+$. This is
possible thanks to a classical result by Carath\'eodory \cite[\S
2]{Pom} concerning the boundary extension of conformal maps between
Jordan domains. Moreover, thanks to an as much as classical result
by Kellogg and Warschawski \cite[\S 3]{Pom}, the map $\varphi_j$
(and its inverse) is of class $C^{2,\alpha}$ up to the boundary.
More precisely, $\varphi_j$ extends to a
$C^{2,\alpha}$-diffeomorphism between $V_j^+\cup\gamma_j$ and
${\mathbb D}^+\cup (-1,1)$. In the following we shall need only the
$C^2$ (or even $C^{1,1}$) regularity, but generally speaking
Kellogg-Warschawski's theorem does not hold in the limit case
$\alpha =0$, whence our assumption $M\in C^{2,\alpha}$ instead of
$M\in C^2$.

Then we define
$$f_j = \varphi_j\circ\pi_j : U_j^+\longrightarrow {\mathbb D}^+$$
where $\pi_j :U_j\to {\mathbb D}$ is the natural projection along
the leaves. Hence $f_j$ is holomorphic in $U_j^+$ and of class $C^2$
up to $M_j\subset \partial U_j^+$. Its differential $df_j$ is a
holomorphic 1-form on $U_j^+$ vanishing on ${\mathcal F}$, thus a
holomorphic section of $N_{\mathcal F}^*$ over $U_j^+$. As a section
of $N_{\mathcal F}^*$, $df_j$ is nowhere vanishing in $U_j^+$.
Moreover, it extends to $M_j$ as a section of class $C^1$, and also
on $M_j$ it is nowhere vanishing.

Finally we define
$$h_j = \log \big\{ \frac{\| df_j\|}{\Im m f_j}\big\} :
U_j^+\longrightarrow {\mathbb R}$$ where $\| df_j\|$ is computed
using the dual norm on $N_{\mathcal F}^*$, induced by the norm with
positive curvature on $N_{\mathcal F}$.

Let us list some properties of these functions $h_j$.

{\bf (1)} First of all, $h_j$ is well defined because $df_j$ does
not vanish on $U_j^+$. Moreover, $\Im m f_j$ tends to zero when
approaching $M_j$, whereas $\| df_j\|$ has a finite nonzero limit,
hence
$$h_j(p) \rightarrow +\infty \qquad {\rm as}\quad p\rightarrow
M_j.$$

{\bf (2)} We may write $h_j = \log \| df_j\| - \log (\Im m f_j)$ and
observe that the second term is plurisubharmonic whereas the ${\rm
i}\partial\bar\partial$ of the first term equals the curvature
$\omega$ of $N_{\mathcal F}$, therefore
$${\rm i}\partial\bar\partial h_j \ge \omega .$$

Now we want to glue the functions $h_j$, and so we need to estimate
their differences $h_j-h_k$ on $U_j^+\cap U_k^+$.

Set $V_{jk}^+ = f_k(U_j^+\cap U_k^+)\subset {\mathbb D}^+$. Then,
for every $j,k=1,\ldots ,\ell$, we have a biholomorphism
$$\varphi_{jk} : V_{jk}^+ \longrightarrow V_{kj}^+$$
such that
$$f_j = \varphi_{jk}\circ f_k \qquad {\rm on}\quad U_j^+\cap
U_k^+.$$

Hence $df_j = (\varphi_{jk}^\prime\circ f_k)\cdot df_k$ and $\Im m
f_j = (\Im m \varphi_{jk})\circ f_k = (\frac{\Im m \varphi_{jk}}{\Im
m}\circ f_k)\cdot \Im m f_k$, and so
$$h_j - h_k = \log \big\{ \vert\varphi_{jk}^\prime\vert\cdot
\frac{\Im m}{\Im m \varphi_{jk}}\big\} \circ f_k .$$

\begin{figure}
\includegraphics[width=8cm,height=5cm]{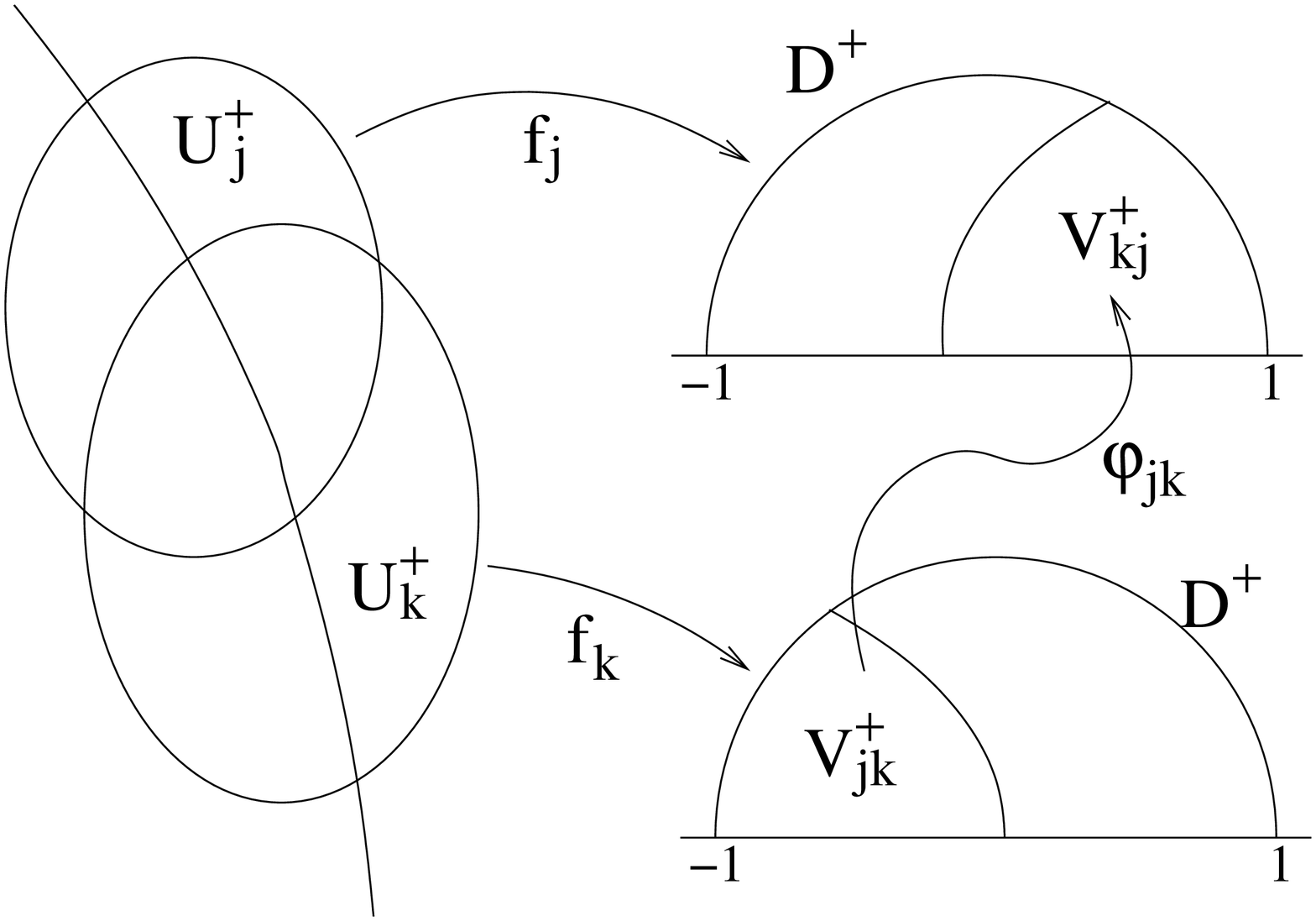}
\end{figure}

Now, by Schwarz reflection the biholomorphism $\varphi_{jk}$ extends
beyond the parts of the boundaries of $V_{jk}^+$ and $V_{kj}^+$
along $(-1,1)$. Moreover, this extension sends the real axis into
itself, and so $\varphi_{jk}^\prime$ is real on $(-1,1)$ and equal
to $\frac{\Im m \varphi_{jk}}{\Im m}$ there. It follows that the
function $\log \big\{ \vert\varphi_{jk}^\prime\vert\cdot \frac{\Im
m}{\Im m \varphi_{jk}}\big\}$ appearing in the above expression of
$h_j - h_k$ is real analytic up to $\partial V_{jk}^+\cap (-1,1)$
and equal to 0 there. As a consequence of this:

{\bf (3)} $$h_j(p)-h_k(p)\rightarrow 0\qquad {\rm as}\quad
p\rightarrow M_j\cap M_k.$$

{\bf (4)} $$dh_j(p) - dh_k(p) \quad {\rm stays\ bounded\ as}\quad
p\rightarrow M_j\cap M_k.$$ Of course, all these estimates hold in a
uniform way. Note that in (4) we use the fact that $f_k$ is of class
$C^1$ up to $M_k$.

We are now ready to construct our exhaustion of $X\setminus M$.

Take a partition of unity $\{ \psi_j\}_{j=1}^\ell$ adapted to $\{
U_j\}_{j=1}^\ell$: each $\psi_j$ is a nonnegative smooth function
with compact support in $U_j$, and $\sum_{j=1}^\ell \psi_j \equiv 1$
around $M$. Define
$$h = \sum_{j=1}^\ell \psi_j h_j : U^+\longrightarrow {\mathbb R}.$$
We obviously have $h(p)\rightarrow +\infty$ as $p\rightarrow M$
(uniformly), by property (1) above.

\begin{lemma} \label{lemma1}
On a sufficiently small neighbourhood of $M$, the function $h$ has
no critical point, and the Kernel $\ker (\bar\partial h)_p\subset
T_pX$ of $(\bar\partial h)_p$ uniformly converges to $T_q^{\mathbb
C}M = T_q{\mathcal F}$ as $p\rightarrow q\in M$.
\end{lemma}

\begin{proof}
Let us work, to fix notation, in the chart $U_\ell$. The
$\bar\partial$-derivative of $h$ can be written as
$$\bar\partial h = \sum_{j=1}^\ell \psi_j\bar\partial h_j +
\sum_{j=1}^{\ell -1} (h_j-h_\ell )\bar\partial\psi_j .$$ The second
term tends to 0 as $p$ tends to $M$, by property (3). For the first
term, we compute
$$\bar\partial h_j = \frac{\bar\partial \| df_j\|}{\| df_j\|} -
\frac{\bar\partial (\Im m f_j)}{(\Im m f_j)} .$$ The 1-form
$(\frac{\bar\partial \| df_j\|}{\| df_j\|})_p$ stays bounded as $p$
tends to $M$, because $f_j$ is of class $C^2$ up to $M$ and $df_j$
does not vanish on $M$. The 1-form $(\frac{\bar\partial (\Im m
f_j)}{(\Im m f_j)})_p$ is, on the contrary, divergent as
$p\rightarrow M$, and moreover its Kernel coincides with the one of
$df_j$, i.e. $T_p{\mathcal F}$. It follows that $\ker (\bar\partial
h)_p$ becomes closer and closer to $T_p{\mathcal F}$ as $p$
approaches to $M$.
\end{proof}

Now we compute the Levi form of $h$ (still in the chart $U_\ell$):
$${\rm i}\partial\bar\partial h = \sum_{j=1}^\ell \psi_j
{\rm i}\partial\bar\partial h_j + \sum_{j=1}^{\ell -1} (h_j-h_\ell
){\rm i}\partial\bar\partial \psi_j + \sum_{j=1}^{\ell -1} {\rm
i}\partial (h_j - h_\ell )\wedge\bar\partial\psi_j +
\sum_{j=1}^{\ell -1} {\rm i}\partial\psi_j\wedge \bar\partial (h_j -
h_\ell )= $$ $$= A + B + C + \overline C .$$ By property (2) we have
$A\ge\omega$, and by property (3) we have $B_p\rightarrow 0$ as
$p\rightarrow M$. By property (4), the 2-form $C_p$ is bounded as
$p\rightarrow M$. Moreover, $C_p$ vanishes on $T_p{\mathcal F}$,
because $\partial (h_j - h_\ell )$ is proportional to $df_\ell$ (see
the computation above). By this and by the previous Lemma,
$C_p\vert_{\ker (\bar\partial h)_p}$ tends to 0 as $p\rightarrow M$.

Therefore, on a sufficiently small neighbourhood of $M$ we certainly
have $${\rm i}\partial\bar\partial h\vert_{\ker (\bar\partial h)} >
\frac{1}{2} \omega\vert_{\ker (\bar\partial h)} .$$ In other words,
the (smooth) hypersurfaces $\{ h=c\}$, $c\gg 0$, are strictly
pseudoconvex. It is then easy to find a convex increasing $r$ such
that $r\circ h$ is strictly plurisubharmonic. After doing the
analogous construction on the negative side $U^-$, we obtain our
desired exhaustion of $X\setminus M$, strictly plurisubharmonic
outside a compact subset.

\begin{remark} {\rm The K\"ahler assumption has not been used up to
now, nor the dimensional assumption $n\ge 3$.}
\end{remark}

\section{End of proof}

Once we know that $X\setminus M$ is strongly pseudoconvex, and $\dim
X\ge 3$, the proof can be concluded in several ways \cite{Lin}
\cite{Siu}.

First of all, we observe that $N_{\mathcal F}\vert_M$ is {\it
topologically} trivial, because it has a nonvanishing section given
by the ``unit normal to $M$'' (as before, we may assume that $M$ is
orientable). Hence, the closed (1,1)-form $\omega$ is exact on $M$,
as well as on a small tubular neighbourhood $U$ of it; we may assume
that $\partial U$ is strictly pseudoconvex (from the exterior
$X\setminus\overline U$), by Proposition \ref{prop1}.

Thus
$$\omega\vert_U = d\beta = \partial\beta^{0,1} +
\bar\partial\beta^{1,0}$$ where the primitive $\beta = \beta^{0,1} +
\beta^{1,0}\in A^1(U)$ can be chosen of real type ($\beta^{1,0} =
\overline{\beta^{0,1}}$) and, from $d\beta = (d\beta )^{1,1}$,
$$\bar\partial\beta^{0,1} = 0.$$

According to a theorem of Grauert and Riemenschneider \cite[Th.
5.8]{Pet}, and because $n=\dim X\ge 3$, the second cohomology group
of the strongly pseudoconvex manifold $X\setminus M$ with ${\mathcal
O}$-coefficients and with compact support $H^2_{\rm cpt}(X\setminus
M,{\mathcal O})$ is equal to zero (indeed, by Serre's duality this
group is isomorphic to the trivial $H^{n-2}(X\setminus M, K_X)$).
This means that the $\bar\partial$-closed $(0,1)$-form
$\beta^{0,1}$, defined on $U$, can be extended to the full $X$, as a
$\bar\partial$-closed $(0,1)$-form $\tilde\beta^{0,1}$: firstly we
extend in any (non $\bar\partial$-closed) way, and then we correct
the error using the vanishing of cohomology with compact support.

Because $X$ is K\"ahler, so that $H^1(X,{\mathcal O})\simeq
H^0(X,\Omega^1)$ by complex conjugation, we may decompose
$$\tilde\beta^{0,1} = \overline\eta + \bar\partial\Phi$$
with $\eta\in\Omega^1(X)$ and $\Phi\in C^\infty (X)$. Hence
$\partial\tilde\beta^{0,1} = \partial\bar\partial\Phi$ and
therefore, setting $\Psi = {\rm i}(\overline\Phi - \Phi )$:
$$\omega\vert_{U'} = {\rm i}\partial\bar\partial\Psi .$$

Thus, we have found a strictly plurisubharmonic function on a
neighbourhood of $M$. But this gives a contradiction with the
maximum principle: a maximum point $p$ for $\Psi\vert_M$ is also a
maximum point for $\Psi\vert_L$, where $L$ is the leaf through $p$,
and this cannot exist.

This end of proof is very close to \cite{Siu}. Really, all the
difficulty of \cite{Siu} is in the fact that there the form $\omega$
is defined only on $M$, and a priori it is not clear how to extend
$\omega$, as a closed (1,1)-form, to a neighbourhood of $M$. Hence
in Siu's paper the (0,2)-form $\bar\partial\beta^{0,1}$ is not
identically zero, as in our case, but it is only vanishing along $M$
at some order (depending on the regularity of $M$). Thus, whereas we
used basically only the resolubility of the $\bar\partial$-equation
with compact support to pass from $\beta^{0,1}$ to
$\tilde\beta^{0,1}$, Siu needs a more delicate result (proved by
himself) on the resolubility of the $\bar\partial$-equation with
growth conditions (see also \cite{Bri}). In our case, as well as in
\cite{Lin}, these difficulties disappear because, by assumption, the
Levi foliation on $M$ can be holomorphically extended to a
neighbourhood of $M$, and this provides the required extension of
$\omega$.

Let us return a moment to the smooth case mentioned in the
Introduction, in absence of a holomorphic foliation. As in
\cite{Siu} and \cite{Bri}, the nonexistence problem is reduced to
construct a plurisubharmonic exhaustion of $X\setminus M$ with some
additional good properties: this permits to prove a suitable
vanishing theorem and then to repeat the arguments above (or,
alternatively, to extend holomorphically the Levi foliation). The
plurisubharmonic exhaustion that we constructed in the previous
Section fits into this general scheme.

A slightly different end of proof is the following one, closer to
the ``topological'' arguments of \cite{Lin}.

Consider the exceptional subset $Y$ of the strongly pseudoconvex
manifold $X\setminus M$. We may find an exhaustion $\psi :
X\setminus M \rightarrow {\mathbb R}$ which is strictly
plurisubharmonic outside $Y$. The classical Morse-type argument of
Andreotti-Frankel-Thom allows to push any compact real surface in
$X\setminus Y$ to a neighbourhood of $M$, using the gradient flow of
$\psi$ (and $n\ge 3$). In other words, $H^2(M,{\mathbb R})$ is
isomorphic to $H^2(X\setminus Y,{\mathbb R})$.

Hence the closed (1,1)-form $\omega$ (which, as before, can be
extended to the full $X$, by pseudoconvexity) is exact not only on
$U$ but even on $X\setminus Y$. By the $\partial\bar\partial$-lemma,
we therefore obtain
$$\omega = \sum_{j=1}^m \lambda_j\delta_{Y_j} + {\rm
i}\partial\bar\partial T$$ where $\{ Y_j\}_{j=1}^m$ are the
irreducible components of $Y$ of codimension one, $\lambda_j$ are
complex numbers, and $T$ is a suitable current, smooth outside
$\cup_{j=1}^m Y_j$. In particular, around $M$ we have $\omega ={\rm
i}\partial\bar\partial T$, and we conclude as before by the maximum
principle.

In fact, this second proof is equivalent to the first one: we have
simply replaced the Hodge symmetry by the
$\partial\bar\partial$-lemma, but the former is also a consequence
of the latter.

\section{Some more remarks}

In trying to extend the previous proof of Theorem \ref{thm1} to the
more general context of Conjecture \ref{conj1}, one is faced with
two main difficulties.

Suppose that Conjecture \ref{conj1} does not hold, and so let
${\mathcal M}\subset X$ be a compact subset invariant by ${\mathcal
F}$ and disjoint from $Sing({\mathcal F})$. We would like to prove
that, thanks to the ampleness of $N_{\mathcal F}$,
$X\setminus{\mathcal M}$ is strongly pseudoconvex.

We cover ${\mathcal M}$ by charts $U_j$, where ${\mathcal F}$ is
defined by $f_j :U_j\rightarrow V_j\subset{\mathbb C}$. Then, on
each $U_j\setminus{\mathcal M}_j$, with ${\mathcal M}_j = {\mathcal
M}\cap U_j$, we take the function
$$h_j = \log \big\{ \frac{\| df_j\|}{dist_j(\cdot ,{\mathcal M}_j)}
\big\}$$ where $dist_j(\cdot ,{\mathcal M}_j)$ is the ``transverse''
distance from ${\mathcal M}_j$, measured with $f_j$, that is
$$dist_j(p,{\mathcal M}_j) = \inf_{q\in{\mathcal M}_j}\vert
f_j(p)-f_j(q)\vert .$$ The functions $h_j$ that we used in Section 2
should be understood as special regularisations of these functions
$h_j$.

It is easily checked that these $\{ h_j\}$ satisfy properties
similar to (1), (2) and (3) of Section 2. However, property (4) is a
more delicate matter, due to the irregular behaviour of
$dist_j(\cdot ,{\mathcal M}_j)$. Let us see an example.

\begin{example} {\rm Take, in the disc ${\mathbb D}$, the closed subset
$K=\{ \arg z =0\ {\rm or}\ \arg z =\frac{\pi}{2} \} \cup \{ 0\}$.
Take two holomorphic embeddings $f_1,f_2 : {\mathbb D}\rightarrow
{\mathbb C}$, and let $g_1,g_2 : {\mathbb D}\rightarrow {\mathbb R}$
be the corresponding distance functions from $K$. Each $g_j$ is not
$C^1$ along an arc $\gamma_j\subset{\mathbb D}$ starting at $0$ with
a tangent of argument $\frac{\pi}{4}$, the equidistant arc from the
two branches of $K$. Typically, these arcs $\gamma_1$ and $\gamma_2$
bound (near the origin) a sector $\Omega$ adherent to $0$, over
which the logarithmic differentials of $g_1$ and $g_2$ are very far
each other: one of them is close to $\frac{dx}{x}$, the other is
close to $\frac{dy}{y}$. Thus $d\log g_1 - d\log g_2$ is unbounded,
on any neighbourhood of $0$.

If we replace $K$ by a curve $K'$ of class $C^{2,\alpha}$, the
situation is not much better: $g_j^2$ are then of class
$C^{2,\alpha}$ up to $K'$, but their quotient is probably no more
than $C^\alpha$ along $K'$, and we don't see how to bound $d\log g_1
- d\log g_2$.}
\end{example}

Related problems appear in trying to extend Lemma \ref{lemma1}.
However, our glueing procedure is rather rudimentary, and one could
suspect that a more refined glueing procedure would work under the
sole assumptions (1), (2), (3) of Section 2.

Suppose now that, in some way, we have proved that
$X\setminus{\mathcal M}$ is strongly pseudoconvex. The second
difficulty is that in Section 3 we used the fact that the 2-form
$\omega$, representing $N_{\mathcal F}$, is exact around $M$. We
don't know if such a fact holds for a more general ${\mathcal M}$.
But \cite{Lin} suggests an alternative approach.

By Baum-Bott formula \cite{Suw}, the cohomology class
$c_1^2(N_{\mathcal F})$ is localized in $Z=Sing({\mathcal F})$. More
precisely, if $\{ Z_j\}_{j=1}^k$ are the irreducible components of
$Z$ of codimension two, then $c_1^2(N_{\mathcal F})$ is cohomologous
to $\sum_{j=1}^k \mu_j[Z_j]$, for suitable complex numbers $\mu_j$
(the Baum-Bott residus along $Z_j$). By the
$\partial\bar\partial$-lemma we therefore have
$$\omega\wedge\omega = \sum_{j=1}^k \mu_j\delta_{Z_j} + {\rm i}
\partial\bar\partial S$$
for a suitable current $S$ of bidegree (1,1).

If $n\ge 3$, the components $Z_j$ are positive dimensional, and
being disjoint from ${\mathcal M}$ they are necessarily contained in
the exceptional subset $Y\subset X\setminus {\mathcal M}$. Hence,
under the canonical contraction $\pi : X \rightarrow X_0$, which
collapses each connected component of $Y$ to a point, each $Z_j$ is
also collapsed to a point, whence the direct image by $\pi$ of the
current $\delta_{Z_j}$ is vanishing. That is,
$$\pi_*(\omega\wedge\omega )={\rm i}\partial\bar\partial S_0$$
where $S_0 = \pi_*(S)$.

This seems a quite strange and unlikely situation. Indeed, $\omega$
is a K\"ahler form (here we are assuming $N_{\mathcal F}$ ample on
the full $X$, not only around ${\mathcal M}$), and we find unlikely
that by a modification the strictly positive (2,2)-form
$\omega\wedge\omega$ may become $\partial\bar\partial$-exact. This
is certainly not the case if $X_0$ is, as $X$, projective. However,
generally speaking $X_0$ is only a Moishezon space \cite{Pet}, non
projective, and on Moishezon spaces we may have nontrivial positive
currents which are $\partial\bar\partial$-exact; but usually (in the
examples we know) these currents do not arise from powers of a
K\"ahler form. Note that the current $\pi_*(\omega^{\wedge n})$ is a
strictly positive measure, hence it cannot be
$\partial\bar\partial$-exact. Also, it is easy to see that
$\pi_*(\omega )$ cannot be $\partial\bar\partial$-exact, otherwise
$\omega$ would be cohomologous to a divisor with support in $Y$, an
evident absurdity. However, the non-$\partial\bar\partial$-exactness
of the intermediate powers $\pi_*(\omega^{\wedge k})$, $2\le k\le
n-1$, seems less evident.

This difficulty does not exist in \cite{Lin} because there
$X\setminus{\mathcal M}$ is not only strongly pseudoconvex but even
Stein, and thus there is nothing to contract.

\end{document}